\documentclass[11pt]{article}
\usepackage{amssymb}
\usepackage{latexsym}

\newtheorem{thm}{Theorem}[section]
\newtheorem{prop}[thm]{Proposition}
\newtheorem{lemma}[thm]{Lemma}
\newtheorem{cor}[thm]{Corollary}
\newtheorem{dfn}[thm]{Definition}

\newtheorem{rmk}[thm]{Remark}

\newcommand{\reals}{\mathbb R}

\newcommand{\G}{{\Lambda}}

\def\qed{\rule{2.3mm}{2.3mm}}

\begin{document}

\title{\bf  A Generalization of Poisson-Nijenhuis Structures }

\author{ 
A\"{\i}ssa WADE 
}

\date{}
\maketitle

\begin{abstract} We generalize  Poisson-Nijenhuis
 structures. We prove that on a manifold endowed with a Nijenhuis tensor and a Jacobi structure which are compatible, there is a hierarchy of 
pairwise compatible Jacobi structures.  Furthermore,  we study the homogeneous
 Poisson-Nijenhuis structures and their relations with Jacobi structures. 

\end{abstract}

\section{Introduction}

 {\it Poisson-Nijenhuis manifolds}, first introduced by Magri and
  Morosi in their paper~\cite{mm} and further studied in~\cite{km},
  play a central role in 
the study of integrable systems. A Poisson-Nijenhuis structure on a manifold 
 $M$ is given by a pair $(\pi, J)$ formed by a Poisson tensor $\pi$ and a 
(1,1)-tensor field $J$  whose Nijenhuis torsion vanishes, such that for
 any two differential 1-forms $\alpha$, $\beta$, we have the following
  compatibility conditions:
  $$\pi(\alpha,\ ^tJ \beta)=  \pi(^tJ\alpha, \ \beta) \quad
{\rm and} \quad C(\pi, J)(\alpha, \beta)=0, \leqno ({\cal C})$$

\noindent where $^tJ$ denotes the transpose of $J$ and $C(\pi, J)$ is an 
 $\reals$-bilinear, skew-symmetric operation on the space
 of differential 1-forms that we shall define below. \par
 In~\cite{v}, the author defined the Poisson-Nijenhuis structures in the
 general algebraic framework of Gel'fand and Dorfman. Moreover, 
 Y. Kosmann-Schwarzbach  gave in~\cite{k} a characterization
 of Poisson-Nijenhuis  structures in terms of Lie algebroids.
  Another characterization is given  in [B-M]. 
 The particular case of a Poisson-Nijenhuis manifold
 having a non degenerate Poisson tensor (or symplectic-Nijenhuis manifold) 
 is especially interesting, this case was studied by several authors
  for different purposes and under various names (see e.g.~\cite{Dorfman}). 
 A natural
 question arises: what is the odd-dimensional analogue of a 
 symplectic-Nijenhuis manifold?  The main aim of the paper 
 is to study this question. A contact manifold is known
 to be the analogue of a symplectic manifold for the odd-dimensional case.
 But a natural framework for a unified study of both contact and
  symplectic manifolds is given by the  Jacobi structures.
 A {\it Jacobi structure} on a manifold $M$ is defined  by a pair $(\G, E)$, where $\G$ is a bivector field, $E$ is a vector field such that $[E,\G]=0$, and 
$[\G,\G]=2E \wedge \G$. Jacobi structures were introduced by A. Lichnerowicz
  and studied by him and his collaborators~\cite{l},~\cite{dlm},~\cite{gl}
  (see also~\cite{Kirillov}).  \

\smallskip
 In the theory of Hamiltonian systems an important mechanism
  allows to construct a hierarchy of pairwise compatible Poisson tensors
  starting from a Poisson-Nijenhuis structure.
 This brings us to the problem  of finding a more general mechanism 
 that would be operational for Jacobi structures. 
 So, we shall need to extend the compatibility 
 conditions $({\cal C})$ above to  the case of Jacobi manifolds. 
 Recently, Marrero, Monterde and Padron (see~\cite{mmp})  considered 
{\it Jacobi-Nijenhuis structures} and gave a possible solution to the
 above questions. Here, we present an approach that is slightly
  different from the one used in~\cite{mmp}. 
 However, these approaches are not equivalent and we shall compare 
 the two approaches in Section 3.2.\par
 
\smallskip
The paper is organized as follows. In Section 2,
 we recall some definitions and basic results concerning Jacobi structures.\par

In Section 3, we  give necessary and sufficient conditions for a
 (1,1)-tensor field $J$ and a Jacobi structure $(\G,E)$ to define, 
 in a natural way,
 a new Jacobi structure which is compatible with $(\G,E)$ in the sense 
 of~\cite{nunes}. This section contains our main results which
  are Theorem~\ref{rec} and Theorem~\ref{pencil}.\par
 
Section 4 is devoted to the analysis of homogeneous Poisson structures,
which are compatible with a (1,1)-tensor field $J$. Such structures are called 
 homogeneous Poisson-Nijenhuis structures. It is well known that homogeneous 
 Poisson structures are related to Jacobi structures, their relations having 
 already  been established in~\cite{dlm}. We give sufficient conditions for the
 existence of a homogeneous Poisson-Nijenhuis structure on a manifold
  and deduce consequences for the existence and properties of
 Jacobi structures.

\section { Preliminaries}
In this paper, all manifolds, multi-vector fields and forms are assumed to be differentiable of class $C^{\infty}$.

\subsection{Jacobi structures}

\begin{dfn} {\rm A {\it Jacobi manifold} $(M,\{ \ , \ \})$ is a manifold $M$ equipped with an $\reals$-bilinear,  skew-symmetric map
 $\{ \ , \ \}:C^{\infty}(M,\reals) \times C^{\infty}(M,\reals) \rightarrow
 C^{\infty}(M, \reals)$, called the {\it Jacobi bracket}, which satisfies the following properties:
\begin{itemize}
\item[1)]  the Jacobi identity: 
$$\{f,\{g,h\}\}+\{g,\{h,f\}\}+\{h,\{f,g\}\}=0, \quad  \forall \ f, g, h \in C^{\infty}(M, \reals);$$
\item[2)] the bracket is local (i.e. the support of $\{f,g\}$ is a subset of the intersection of the supports of $f$ and $g$).
\end{itemize}  
}
\end{dfn}
 Equivalently, a Jacobi structure can be defined by  a pair $(\G, E)$ 
 of a bivector field $\G$ and a vector field $E$ such that 

$$[E,\G]=0 \quad {\rm and} \quad [\G, \G]=2E \wedge \G,$$ 

\noindent where $[ \ , \ ]$ is the Schouten-Nijenhuis bracket on the space of
multi-vector fields (see~\cite{kz}). The Jacobi bracket is then given by
$$\{f,g\}=\G(df,dg)+\langle fdg-gdf, E\rangle.$$ 

 When $E$ is zero, we obtain a Poisson structure. In other words, a {\it Poisson structure} on a manifold $M$ is given by a bivector field $\G$ such that the Schouten-Nijenhuis bracket $[\G,\G]$ vanishes. Then $(M,\G)$  is called a
 {\it Poisson manifold}. In~\cite{l}, Lichnerowicz  showed that to any Jacobi structure
$(\G,E)$ on a manifold $M$, one may associate a Poisson structure
 $\pi$ on $M\times\reals$, defined by

$$\pi(x,t)=e^{-t}\Big( \G(x)+{\partial \over \partial t} \wedge E \Big).$$

\noindent Then, $\pi$ is called the {\it Poissonization} of $(\G, E)$. 

\bigskip

 \noindent {\bf Notations.} To any bivector field $\G$ on $M$, we may associate the skew-symmetric linear map denoted also by $\G :\Omega^1(M) \rightarrow \chi(M)$ and defined by:

$$\langle \beta, \ \G \alpha \rangle= \langle \alpha \wedge \beta, \ \G\rangle=
 \G(\alpha, \beta),$$

 \noindent where $\chi(M)$ is the space of vector fields
and $\Omega^1(M)$ is the space of differential 1-forms on $M$.
 Conversely, a linear map $\G :\Omega^1(M) \rightarrow \chi(M)$ defines a bivector field on $M$ if and only if 
$$\langle \alpha , \G\beta \rangle + \langle \beta, \G \alpha \rangle=0.$$

\bigskip

\noindent {\bf Example: locally conformal symplectic manifolds.} Let $M$ be a $2n$-dimensional manifold. A {\it locally conformal symplectic structure} 
on $M$ is given by a pair $(F, \omega)$, where $F$ is a non-degenerate 2-form and $\omega$ is a 1-form 
such that  

$$ d \omega=0 \quad {\rm and} \quad dF+ \omega \wedge F=0.$$

\noindent We define a bivector field $\G$ and a vector field $E$ by:

$$ i_E F= \omega \quad {\rm and} \quad i_{\G \alpha}F= - \alpha, \ \ \forall 
 \alpha \in \Omega^1(M).$$

\noindent Then $(\G,E)$ defines a Jacobi structure. For any 
$x  \in M$, there exist a neighborhood $U_x$ and a function $f$ defined on 
$U_x$ such that $\omega=df$ and $\Omega=e^f F$ is symplectic. 

\bigskip

\subsection{The Lie algebroid of a Jacobi manifold}

It was proven in [Ke-SB] that there is a Lie algebroid associated 
 with an arbitrary Jacobi manifold $(M,\G, E)$.
 We refer the reader to~\cite{Mackenzie}, for instance, for the basic
 properties of Lie 
 algebroids.  Consider the vector bundle $T^*M \oplus \reals$. The space $\Gamma(T^*M \oplus \reals)$ of smooth sections may be identified with $ \Omega^1(M) \times 
C^{\infty}(M , \reals)$.  
The Lie algebroid associated with a Jacobi manifold $(M,\G, E)$ is $T^*M \oplus \reals$ with  the Lie bracket $\{ \ , \ \}_{(\G,E)}$ on 
$\Gamma(T^*M \oplus \reals)$ defined by

\begin{eqnarray*}
\{(\alpha,f),(\beta ,g)\}_{(\G, E)}&=\Big(L_{\G \alpha}\beta - L_{\G \beta}\alpha -d \big(\G(\alpha, \beta)\big) +fL_E\beta -gL_E\alpha -i_E(\alpha \wedge \beta), \\
&  -\G(\alpha, \beta) + \G(\alpha, dg)-\G(\beta, df) +fE(dg) -gE(df) \Big),
\end{eqnarray*}

\noindent where $d$ is the exterior derivative and $L_X=di_X + i_Xd$ is the Lie derivation by $X$, for any vector field $X$. 
The anchor is given by the map $\#_{_{(\G,E)}}$ such that 

$$\#_{_{(\G,E)}}(\alpha,f)= \G \alpha+fE.$$

\section {Jacobi and Nijenhuis structures}
\subsection{Extension of the definition of compatibility  to Jacobi structures}

 Let $J$ be a $(1,1)$-tensor field of $M$.
The {\it Nijenhuis torsion} $N_J$ of $J$ with
respect to the Lie bracket $[.,.]$ on the space $\chi(M)$ of vector fields is defined by
$$N_J(X,Y)=[JX, JY] - J[JX, Y] - J[X, JY]+J^2[X,Y], \quad \forall \ X, Y \in \chi(M).$$

\begin{dfn}  $J$ is called a {\it Nijenhuis tensor} if its Nijenhuis torsion vanishes. 
\end{dfn}

 In particular, when  $J$ is a $(1,1)$-tensor field on $M$ and 
 $\G : \Omega^1(M) \rightarrow \chi(M)$ is a linear map, then $J \circ \G$
defines a bivector field if and only if $J\circ \G= \G\circ \ ^tJ$.
In this case, the associated bivector field is denoted by $J\G$. \par

\smallskip

 Furthermore, any bivector field $\G$ gives rise to a bracket defined on the 
 differential 1-forms by 

\begin{equation} \label{eq:bra}
\{\alpha, \beta\}_{\G}=L_{\G \alpha}\beta -L_{\G \beta} \alpha -d \big(\G(\alpha, \beta)\big), \quad \forall \alpha, \beta \in \Omega^1(M),
\end{equation}

\noindent where $L_X$ is the Lie derivation by $X$, for any vector field $X$.

\smallskip

 Whenever $J\circ \G=\G \circ \ ^tJ$, we denote by $C(\G, J)$ the $\reals$-bilinear map given  by 
$$C(\G,J)(\alpha,\beta)= \{\alpha,\beta\}_{J\G}-
\Big(\{ \ ^tJ\alpha, \beta\}_{\G}+\{\alpha,\ ^tJ\beta\}_{\G}-
 \ ^tJ\{\alpha, \beta\}_{\G}\Big).$$ 

\begin{dfn}(see~\cite{km}) A  Poisson-Nijenhuis structure on a manifold 
$M$ is defined by a Poisson tensor $\pi$ and a Nijenhuis tensor $J$ on $M$ 
 such that 
\begin{itemize}
\item[(1)] $J\circ \pi=\pi \circ \ ^tJ,$
\item [(2)] $C(\pi,J)=0$.
\end{itemize}
In this case, we say that $\pi$  and $J$ are compatible.
\end{dfn}

 An extension of this definition to Jacobi structures is given by:

\begin{dfn}
\label{Com} Let $(M,\Lambda, E)$  be a Jacobi manifold and let
 $J$ be a $(1,1)$-tensor field such that its Nijenhuis torsion $N_J$ satisfies 
$$N_J(\G \alpha,\G \beta)=0 \quad {\rm and} \quad
N_J(\G \alpha, J^kE)=0,$$

\noindent for all $\alpha, \  \beta \in \Omega^1(M)$ 
  and for all $k\in \mathbb N$. 
 Then, $J$ is said to be compatible with the Jacobi 
 structure $(\Lambda,E)$ if 
\begin{itemize}
\item[(i)]  $J\circ \G=\G \circ \ ^tJ,$ 
\item [(ii)] $ \big \langle \ ^t J \gamma, \ \G \Big(C(\G,J)(\alpha, \beta) 
 \Big) - \G(\alpha, \beta)JE +\G(\alpha, \ ^t J\beta)E \big \rangle=0,$
for any $\alpha, \  \beta, \  \gamma \in \Omega^1(M)$,
\item[(iii)] $ [J^kE, \G ]+[E, J^k\G]=0,$ for any integer $k \geq 1$.
\end{itemize}
\end{dfn}

 When $E=0$ (i.e. $\G$ defines a Poisson structure), the definition 
 reduces to that of a {\it weak Poisson-Nijenhuis structure} (see~\cite{mn}). 
 The compatibility conditions then reduce  to 
 $(i)$ and $(ii)$ and constitute a generalization of the definition of a 
 Poisson-Nijenhuis structure introduced in ~\cite{mm}.

\begin{thm} 
\label{rec}
Let $(\G,E)$ be a Jacobi structure  on $M$. Assume that
$J$ is a $(1,1)$-tensor field  such that $J\circ \G=\G \circ \ ^tJ$,  

$$N_J(\G \alpha,\G \beta)=0
\quad {\rm and} \quad
N_J(\G \alpha, E)=0 \quad \forall \alpha, \beta \in \Omega^1(M),$$
 \noindent where $N_J$ is the Nijenhuis torsion of $J$. Then $(J \G, JE)$ is a Jacobi structure on $M$  if and only if the following properties are satisfied
for all $\alpha, \beta, \gamma \in \Omega^1(M):$

\begin{itemize}
\item[(a)] $ J([JE, \G ]\alpha+[E, J\G]\alpha)=0,$
\item [(b)] $ \big \langle \ ^t J \gamma, \ \G \Big(C(\G,J)(\alpha, \beta) \Big) - \G(\alpha, \beta)JE +\G(\alpha, \ ^t J\beta)E \big \rangle=0.$ 
\end{itemize}

In particular, if $J$ is  a Nijenhuis tensor compatible with $(\G,E)$, then $(J \G, JE)$ is a Jacobi structure on $M$.

\end{thm}

The proof of Theorem~\ref{rec} is based on the following three lemmas.

\begin {lemma}
\label{no-pr}
For any bivector field $\G$, 
\begin{equation}
 \big \langle \gamma, \G \{\alpha ,\beta\}_{\G} \big \rangle=
 \big \langle \gamma, [\G \alpha ,\G \beta] \big \rangle+ {1 \over 2}
[\G,\G](\alpha, \beta, \gamma), \quad 
\forall  \ \alpha, \beta, \gamma \in \Omega^1(M).
\end{equation}

\end{lemma}
\noindent This formula is proven for instance in [K-M]. 

\begin{lemma}
\label{bivec}
Consider a pair $(\G , E)$ of a bivector field $\G$ and a vector field
 $E$ on $M$ such that $[\G,\G]= 2E \wedge \G$. Then, for any 
 $C^{\infty}(M,\reals)$-linear map $J$ on
 $\chi(M)$ satisfying $J \circ \G= \G \circ \ ^t J$, the following formula holds:

\begin{eqnarray*}
{1 \over 2}[J\G, J\G](\alpha, \beta, \gamma)& = &  \ ( JE \wedge  J\G) (\alpha, \beta, \gamma)+ \big \langle ^tJ \gamma, \ \G \Big( C(\G, J)(\alpha, \beta)\Big)  \big \rangle  \\ 
& &+E(^tJ \gamma) \G(\alpha, \ ^tJ \beta)-JE(^tJ \gamma) \G (\alpha, \beta)  \\ & & - \big \langle \gamma, \ N_J(\G \alpha, \G \beta)  \big \rangle.
\end{eqnarray*}

\end{lemma}

\noindent {\it Proof:}  Lemma~\ref{no-pr} yields

$${1 \over 2} [J\G  ,J\G](\alpha, \beta, \gamma)=
 \big \langle ^tJ \gamma, \ \G \{\alpha ,\beta\}_{J\G}\big \rangle - \big 
\langle \gamma, \ [J\G \alpha ,J\G \beta] \big \rangle.$$

\noindent A direct computation, using  Relation (2) again, shows that

\begin{eqnarray}
{1 \over 2} [J\G  ,J\G](\alpha, \beta, \gamma)&=&
{1 \over 2} \bigg([\G, \G](^tJ\alpha,\beta, \ ^tJ \gamma) 
+[\G, \G](\alpha, \ ^tJ\beta, \ ^t J \gamma) \nonumber \\
  & &-[\G, \G](\alpha, \beta ,\  ^tJ^2 \gamma)\bigg) 
 + \langle ^tJ \gamma, \G \Big(C(\G, J)(\alpha, \beta) \Big)\rangle  \nonumber  \\
& & - \big \langle  \gamma, N_J(\G \alpha, \G \beta) \big \rangle.
\end{eqnarray}

\noindent Since $[\G, \G]= 2 E\wedge \G$, we obtain: 
\begin{eqnarray*}
 {1 \over 2} [J\G  ,J\G](\alpha, \beta, \gamma) &=  
 (J E \wedge J  \G) (\alpha, \beta, \gamma) 
+E(^tJ \gamma) \G(\alpha, \ ^tJ \beta)
 -JE(^tJ \gamma)  \ \G (\alpha, \beta) \\
&  + \big \langle^t J \gamma, \ \G \big( C(\G, J)(\alpha, \beta) \big)
\big \rangle - 
 \big \langle  \gamma, N_J(\G \alpha, \G \beta) \big \rangle .
\end{eqnarray*}

\noindent This is the formula that we sought.
\hfill \qed 

\begin{lemma}
\label{vec}
Let $\G$ be a bivector field and $E$ a vector field on $M$.
 The following relation holds for any $C^{\infty}(M,\reals)$-linear 
 map $J$ on $\chi(M)$:
$$[JE, J \G](\alpha, \beta)= \langle \beta, \ N_J(E,\G \alpha) \rangle+
\langle \beta, J[JE, \G]\alpha+ J[E,J\G] \alpha  -J^2[E, \G] \alpha\rangle.$$

\end{lemma}

\noindent {\it Proof:} For any bivector field $\G$ and for all $\alpha , \beta \in
\Omega^1(M)$, we have:
$$[E,\G](\alpha, \beta)=L_E\big(\G(\alpha,\beta) \big)- 
\G(L_E \alpha, \beta)- \G(\alpha, L_E \beta).$$
 
\noindent This is equivalent to the relation  
\begin {equation}\label{eq:formula4}
[E, \G]\alpha=[E,\G \alpha]-\G L_E \alpha,\quad \forall \alpha \in \Omega^1(M). \end{equation} 

\noindent Using (4), we obtain for any $\alpha \in \Omega^1(M)$:

\begin{eqnarray*}
[JE, J\G] \alpha &=& [JE, J\G\alpha] - J\G L_{JE} \alpha \\
& =& N_J(E, \G \alpha)+J[JE, \G \alpha]+ J[E, J\G \alpha] -J^2[E, \G \alpha]
-J \G L_{JE} \alpha.
\end{eqnarray*}

\noindent Replacing $[E, \G \alpha]$ by $ [E, \G]\alpha +\G L_E \alpha$, we deduce that
\begin{eqnarray*}
[JE, J\G] \alpha &= &N_J(E, \G \alpha) +J \big([JE, \G \alpha]-
\G L_{JE} \alpha \big) \\
& &   + J\big([E,J\G \alpha]-J \G L_E \alpha \big) -J^2 [E,\G] \alpha  \\
&=&N_J(E, \G \alpha)+ J\big([JE, \G]+ [E,J\G]-J[E, \G]\big) \alpha.
 \end{eqnarray*}
 
\hfill \qed

\medskip
\noindent {\it Proof of Theorem~\ref{rec}:} Lemma~\ref{vec} ensures that $[JE,J\G]=0$ is equivalent to (a), while
Lemma~\ref{bivec} asserts that $[J\G, J\G]= 2 JE \wedge J\G$  if and only if 
 property (b) is satisfied. Therefore, the theorem is proved. 
\hfill \qed

\bigskip

 Now, let us  express Properties (a) and (b) of Theorem~\ref{rec}
 using the Lie algebroid bracket associated with the Jacobi structure. 

\begin{prop} 
Let  $(\G, E)$ be a Jacobi structure on $M$ and let $J$ be a $(1,1)$-tensor
 field  such that  
$$J\circ \G=\G \circ \ ^tJ\quad {\rm and}\quad 
N_J(\G \alpha+fE,\G \beta+gE)=0, \quad \forall \alpha, \beta \in \Omega^1(M), 
\ \forall f,g \in C^{\infty}(M, \reals).$$

\noindent Then,

\begin{eqnarray*}
 (a) \ \mbox{\rm is satisfied }  &\iff& 
 [J\G\alpha+f JE, \ gJE]=\#_{_{(J\G,JE)}} \Big(\{(\alpha,f),(0,g)\}_{(J\G,JE)}
  \Big)\\
 (b) \ \mbox {\rm is satisfied } &\iff&
  [J\G\alpha, \ J\G \beta]=\#_{_{(J\G,JE)}} \Big(\{(\alpha,0),(\beta,0)\}_{(J\G,JE)}\Big).
   \end{eqnarray*}

\end{prop}

\noindent{\it Proof:} On one hand,
$$[J\G \alpha+f JE, \ gJE ]= g[J \G\alpha, JE]+\big(J \G(\alpha,dg)+
\langle fd g -gdf , JE \rangle \big) JE.$$

\noindent On the other hand, 
$$\#_{_{(J\G,JE)}}\{(\alpha,f), (0,g)\}_{(J\G,JE)}=- gJ \G L_{JE} \alpha 
+\big(J \G(\alpha,dg)+ \langle f dg -gdf , \ JE \rangle \big) JE.$$

\noindent We deduce that

\begin{eqnarray*}
[J\G\alpha+f JE, gJE]-\#_{_{(J\G,JE)}}\{(\alpha,f), (0,g)\}_{(J\G,JE)}&=&
g([J \G\alpha, JE]+J \G L_{JE} \alpha) \\
& =&g [J \G, JE]\alpha.
\end{eqnarray*}

\noindent But Lemma~\ref{vec} asserts that   

$$[J \G, JE]\alpha=0 \iff J([JE, \G]\alpha + [E,J \G] \alpha)=0.$$

\noindent Hence we obtain the first equivalence. In the same way, we prove the
second equivalence using Lemma~\ref{bivec}.

\subsection{ Hierarchy of Jacobi structures}

\indent A manifold $M$ is said to be a {\it bihamiltonian manifold} if $M$ is endowed with two Poisson tensors $\pi_1$ and $\pi_2$ such that 
$\pi_1-\lambda \pi_2$ is a Poisson tensor for any $\lambda \in \reals$. Then
$\pi_1-\lambda \pi_2$ is called a {\it Poisson pencil}.
By analogy, if $\{.,.\}_{_1}$ and  $\{.,.\}_{_2}$  are two Jacobi structures
such that  $\{.,.\}_{_{\lambda}}=\{.,. \}_{_1} - \lambda  \{.,.\}_{_2}$  defines a Jacobi structure for any $\lambda$ in $\reals$, then 
 $\{.,.\}_{_\lambda}$ will be called a {\it Jacobi pencil}. In this case, 
 the two Jacobi structures are said to be {\it compatible} (see [N]).

\begin{prop}
\label{car} Let $(\Lambda_1, E_1)$ and $(\Lambda_2, E_2)$ be two Jacobi 
 structures on $M$. Denote by 
$\pi_i=e^{-t}(\G_i + \partial/ \partial t \wedge E_i)$,  with $i=1,2$, the associated Poisson tensors on $M \times \reals $.
Then the following assertions are equivalent: 

\begin{itemize}
\item[(1)] $(\Lambda_1, E_1)$ and $(\Lambda_2, E_2)$ define a Jacobi pencil on $M$.
\item[(2)] $[E_1, \G_2]+[E_2, \G_1]=0$ and 
 $[\G_1,\G_2]=E_1\wedge \G_2 +E_2\wedge \G_1$.

\item[(3)] The pair $(\pi_1, \pi_2)$ defines a Poisson pencil on $M \times \reals $.

\end{itemize}
\end{prop}

\noindent The proof of this proposition is straightforward and can be found in   [N].
 
\begin{rmk}
\label{remarque}
{\rm Consider the following relation:

\begin{itemize}
 \item[($a'$)] \ \ \ \ \ \ \ \ \ \ \ $[JE, \G]+ [E,J \G]=0$ 
\end{itemize}

\noindent Under the hypotheses of Theorem~\ref{rec}, if $(J\G, JE)$ is a 
 Jacobi structure compatible with $(\G,E)$, then it follows 
 from the proposition 
 above and Theorem~\ref{rec} that ($a'$), as well as  $(b)$ are satisfied.\par 
 On the other hand, if  ($a'$) and $(b)$ are
 satisfied then Theorem~\ref{rec} ensures that $(J\G, JE)$  is a Jacobi
 structure on $M$. But, $(J\G, JE)$ may not be compatible with $(\G,E)$.
 They define a Jacobi pencil if and only if 
 $[J \G , \G]=JE \wedge \G + E \wedge J\G.$

}
\end{rmk}

\begin{thm} 
  \label{pencil} For any Jacobi structure  $(\G, E)$ compatible with a Nijenhuis tensor $J$ on $M$ and  for each $k \in \mathbb N^*$, the pair $(J^k\G, J^k E)$ is a Jacobi structure on $M$. 
Furthermore for $k_1,$ $k_2 \in \mathbb N^*$,  $(J^{k_1}\G,J^{k_1} E)$ and  $(J^{k_2}\G, J^{k_2} E)$ define a Jacobi pencil.
\end{thm}

\noindent This theorem is a generalization of a result proved in 
[M-M] and [K-M].

\begin{lemma}
\label{nij} Let $J$ be a $(1,1)$-tensor field. Then, we have:
\begin{eqnarray*}
 N_{J^{k+1}}(X, Y) &=&  N_{J^k}(JX ,JY)
+ J^k \Big( N_J(J^kX, Y)+ N_J(X ,J^k Y) \Big) \\
 & & -J^2 \Big(  N_{J^{k-1}}(JX, JY) -  N_{J^k}(X , Y) \Big), \quad \forall X,Y \in \chi(M).
\end{eqnarray*}
\end{lemma}

\noindent The proof of this lemma is straightforward.

\medskip
\noindent {\it Proof of Theorem~\ref{pencil}:} assume that
 $[J^{\ell}\G,J^{\ell}\G]=2J^{\ell} E \wedge J^{\ell}\G$, for any $\ell \leq k$.
It follows from Lemma~\ref{bivec} that 
\begin{eqnarray*}
{1 \over 2}[J^{k+1}\G, J^{k+1}\G](\alpha, \beta, \gamma)& =
(J^{k+1}E \wedge  J^{k+1}\G) (\alpha, \beta, \gamma)+ \langle ^tJ \gamma, \ 
J^{k}\G C(J^{k}\G, J)(\alpha, \beta) \rangle  \\ 
& +J^{k}E(^tJ \gamma) J^{k}\G(\alpha, \ ^tJ \beta)-J^{k+1}E(^tJ \gamma) 
J^{k}\G (\alpha, \beta),
\end{eqnarray*}

\noindent for all for all $\alpha, \ \beta, \ \gamma$ in $\Omega^1(M)$.
 We shall prove that 
$$ J^{k}\G C(J^{k}\G, J)(\alpha, \beta)+J^{k}\G(\alpha, \ ^tJ \beta)J^{k}E -
J^{k}\G (\alpha, \beta)J^{k+1}E=0, $$
 
\noindent for any $k \geq 1$. In fact,  for any bivector field $\G$  and  for any linear map $J$ on
 $\chi(M)$  such that $J \circ \G=\G \circ \ ^tJ$, the following relation holds (see [M-M]):
\begin{equation}
\label{cj}
 \langle C(J \G, J)( \alpha, \beta), \ X \rangle= \langle C(\G,J)(^tJ \alpha, \beta), \ X \rangle+ \langle \alpha, \ N_J(\G \beta ,X)\rangle.
\end{equation}

\noindent  For any vector field $X$ of the form $\G \gamma$, we have

$$\langle C(J \G, J)( \alpha, \beta), \ \G \gamma \rangle= \langle C(\G,J)(^tJ \alpha, \beta), \ \G \gamma \rangle$$

\noindent Thus,

$$\G \ C(J\G, J)(\alpha, \beta)=\G \ C(\G, J)(^tJ\alpha, \beta),$$
\noindent  for all $\alpha, \ \beta$ in $\Omega^1(M)$. 
Hence, we obtain by induction that for any 
$k \geq 1$, 

\begin{equation}
\label{cj2}
 J^k\G \ C(J^k\G, J)(\alpha, \beta) =J^k\G \ C(\G, J)(^tJ^k\alpha, \beta).
\end{equation}

\noindent Since $J$ is compatible with $(\G,E)$, we have 

\begin{equation}
\label{cj3}
 J \G \ C(\G,J)(\alpha, \beta) =J\Big(\G (\alpha, \beta)JE -
\G(\alpha, \ ^tJ \beta) E\Big) .
\end{equation}

\noindent We deduce that 

\begin{eqnarray*}
J^{k}\G C(J^{k}\G,J)(\alpha, \beta)&=& J^{k}\G C(\G, J)(^tJ^k\alpha, \beta)\\
 &=& J^{k}\Big(\G(^tJ^k\alpha, \beta)JE - \G(^tJ^k\alpha, ^tJ\beta)E\Big) \\
&=& J^{k}\G (\alpha, \beta)J^{k+1}E-J^{k}\G(\alpha, \ ^tJ \beta)J^{k}E,
\end{eqnarray*}

\noindent for any $k \geq 1.$  So,  we obtain the relation that we sought. The latter implies that 
$$[J^k \G, J^k \G]= 2J^k E \wedge J^k \G, \quad {\rm for \ any} \quad k \geq 1.$$
 
 \noindent Moreover, replacing $J$ by $J^k$ in Lemma~\ref{vec}, we obtain,
 since $[E, \G]=0$, 
$$[J^kE, J^k \G](\alpha, \beta)= \langle \beta, \ N_{J^k}(E,\G \alpha) \rangle+
\langle \ ^t J^k \beta, [J^kE, \G]\alpha+ [E,J^k\G ]\alpha \rangle.$$

\noindent From Lemma~\ref{nij}, we obtain by induction that
 $N_{J^k}(E,\G \alpha)$  vanishes for all $k \geq 1$. Therefore, 
$$[J^kE, J^k \G]=0, \quad  {\rm for \ all} \quad k \geq 1.$$ 
 
\noindent Thus, $(J^k \G, J^kE)$ defines a Jacobi structure for any 
$k \geq1$.

\smallskip
Now, take two different pairs $(J^{k_1}\G,J^{k_1} E)$ and  $(J^{k_2}\G, J^{k_2} E)$. We shall prove that they determine a Jacobi pencil. 
For any $\lambda \in \reals$, we have to prove that
$$[J^{k_1} \G- \lambda J^{k_2} \G, J^{k_1} \G- \lambda J^{k_2} \G]= 2
(J^{k_1}E-\lambda J^{k_2} E) \wedge (J^{k_1} \G-\lambda J^{k_2} \G).$$

\noindent We already know that 
$$[J^{k_i} \G, J^{k_i} \G]= 2J^{k_i} E \wedge J^{k_i} \G, \quad \forall i=1,2.$$

\noindent Now we prove that 
$$[J^{k_1} \G,J^{k_2} \G]= J^{k_1} E \wedge J^{k_2} \G + J^{k_2} E \wedge 
J^{k_1}\G .$$

\noindent Assume that $k_1=k_2+ \ell$, then we apply $\ell$ times
 the result saying that, for arbitrary bivector fields $\G$ and $\pi$ on $M$,
 for any linear map $J$ on $\chi(M)$ the following formula holds (see [M-M]):

\begin{eqnarray}
 \label{eq:mm}
 [J \G, \pi]( \alpha, \beta, \gamma)&= & [\G, \pi]( \alpha, \beta, \ ^tJ\gamma)  +\langle C( \pi, J)( \alpha, \gamma) , \ \G \beta \rangle  \nonumber \\
& & - \langle C( \pi, J)( \beta, \gamma) , \ \G \alpha \rangle 
- \langle C(\G, J)( \alpha, \beta) , \ \pi \gamma \rangle.
\end{eqnarray}
 
\noindent We apply this last relation and we calculate by recursion the $\ell$
quantities $[J^{k_2+ \ell} \G , J^{k_2} \G]$,...,$[J^{k_2+ 1} \G, J^{k_2} \G]$.   It follows that: 

\begin{eqnarray*}
[J^{k_1} \G, J^{k_2} \G]( \alpha, \beta, \gamma) &=&
[J^{k_2} \G, J^{k_2} \G]( \alpha, \beta,\  ^tJ^{\ell}\gamma) \\
 && +\sum_{r=1}^{\ell}\langle C( J^{k_2}\G, J)( \alpha, \ ^tJ^{r-1}\gamma) , \
 J^{k_1 -r}\G \beta \rangle \\
& & - \sum_{r=1}^{\ell}\langle C( J^{k_2}\G, J)( \beta, \ ^tJ^{r-1}\gamma) , \
 J^{k_1 -r}\G \alpha \rangle \\
&&  -  \sum_{r=1}^{l}\langle C( J^{k_1-r}\G, J)(\alpha, \beta)           , \ J^{k_2 +r-1}\G \gamma \rangle                                             \end{eqnarray*}

\noindent Using (~\ref{cj2}),  (~\ref{cj3}) and the fact that
$[J^{k_2} \G, J^{k_2} \G]=2J^{k_2} E \wedge J^{k_2} \G$, we obtain after 
computations:

$$[J^{k_1} \G, J^{k_2} \G] = J^{k_1} E \wedge J^{k_2} \G 
+J^{k_2} E \wedge J^{k_1} \G.$$

\noindent The last step is to show that:
$$[J^{k_1} E- \lambda J^{k_2} E, J^{k_1} \G- \lambda J^{k_2} \G]= 0.$$

\noindent This is equivalent to showing that  $[J^{k_1} E, J^{k_2} \G]+ 
[ J^{k_2} E,J^{k_1} \G]= 0.$ By hypothesis this relation is true when 
$k_2=1$ and using Lemma~\ref{vec}, we can easily show by induction that 
this formula holds for any $k_1$ and $k_2$. 
\hfill \qed

\bigskip
\noindent{\bf Example.} Let  $\omega$ be a closed 1-form and let $F_1$, $F_2$
be two non-degenerate 2-forms on $M$. Assume that $(F_1,\omega)$  and 
$(F_2,\omega)$ are  locally conformal symplectic structures on $M$. 
 Let $(\Lambda_i, E_i)$ denote the Jacobi structures associated with $(F_i,\omega)$, where $i=1,2$. Assume that these two Jacobi structures are compatible. Define the isomorphism  of $C^{\infty}(M, \reals)$-modules $\flat_i : \chi(M) \rightarrow \Omega^1(M)$ by
 $$\flat_i(X)=  -i_XF_i. $$ 

\noindent We have 
$$E_i= -\flat_i^{-1}(\omega) \quad {\rm and } \quad \G_i \alpha= \flat_i^{-1}(\alpha), \quad \forall \alpha \in \Omega^1(M).$$

\noindent Then, the $(1,1)$-tensor field $J= \flat_2^{-1} \circ \flat_1$
 is compatible with $(\Lambda_1, E_1)$. Indeed, for any $x  \in M$, there exist a neighborhood $U_x$ and a function $f$ defined on $U_x$ such that $\omega=df$. The 2-forms  $\Omega_1=e^f F_1$ and  $\Omega_2=e^f F_2$ are 
symplectic and the Poisson tensors associated with $\Omega_1$, $\Omega_2$ are
respectively $\pi_1=e^{-f} \G_1$, $\pi_2=e^{-f} \G_2$.

 We claim that the Jacobi structures $(\Lambda_1, E_1)$  and $(\Lambda_2, E_2)$ are compatible if and only if $\pi_1$ and $\pi_2$ 
are compatible. Let us prove this claim. Using  the properties of the 
Schouten-Nijenhuis bracket, we get 
 
$$
[\pi_1,\pi_2]= e^{-2f}\Big ([\G_1, \G_2] -[\G_1, f] \wedge \G_2
- [\G_2, f] \wedge \G_1 \Big).$$

\noindent Since $E_i=[\G_i,f]= - \G_i (df)$, we have

$$ [\pi_1,\pi_2]= e^{-2f}\big ([\G_1, \G_2] -E_1 \wedge \G_2
- E_2 \wedge \G_1 \big).$$

\noindent Therefore, $[\pi_1,\pi_2]=0$ if and only if 
$[\G_1, \G_2] =E_1 \wedge \G_2 +E_2 \wedge \G_1$. Moreover, we may remark that the Jacobi identity for the Schouten-Nijenhuis 
bracket yields

\begin{eqnarray*}
[[\pi_1,\pi_2], e^f] &=& -[[ \pi_2, e^f], \pi_1] - [[ e^f,\pi_1], \pi_2]\\
&=&-[[ \G_2, f],e^{-f} \G_1] - [[ f,\G_1], e^{-f} \G_2].
\end{eqnarray*}

\noindent The fact that  $E_i=[\G_i,f]$ implies 

$$[[\pi_1,\pi_2], e^f]= -e^{-f} ([E_2, \G_1] +[E_1, \G_2]). $$

\noindent Thus, $(\G_i, E_i)_{i=1,2}$ form a Jacobi pencil 
if and only if the tensors $(\pi_i)_{i=1,2}$ define a Poisson pencil. So, we may deduce that the Nijenhuis torsion of $J$  vanishes. 
 Furthermore, the sequence $(J^k \pi_1)$ consists of pairwise compatible
Poisson tensors, while $( J^k \G_1, J^k E_1)$ is a sequence of pairwise 
compatible Jacobi structures.

\subsection {Compatibility conditions and Jacobi-Nijenhuis structures}

  In this section, we shall study the differences between a Jacobi-Nijenhuis
 structure and  a structure which satisfies the axioms of
  Definition~\ref{Com}. We now recall
  the definition of a Jacobi-Nijenhuis structure given in~\cite{mmp}.
  For any bivector field $\G$ and for any vector field $E$, we may define 
 a map $\widetilde \#_{_{(\G,E)}}:
 \Omega^1(M)\times C^{\infty}(M,\reals)\rightarrow \chi(M) \times 
 C^{\infty}(M, \reals)$ by

$$\widetilde \#_{_{(\G,E)}}(\alpha,g)=(\G \alpha+gE, - i_E\alpha), \quad
 \forall \ \alpha \in \Omega^1(M), \ \forall \ g \in  C^{\infty}(M,\reals).$$

 \noindent  Consider the Lie bracket $[.,.]_{\star}$ defined 
 on  $\chi(M) \times C^{\infty}(M, \reals)$  by:
 $$[(X_1,f_1) , \ (X_2,f_2)]_{\star}= ([X_1,X_2], i_{X_1}df_2- i_{X_2}df_1).$$  
 
\noindent  Let ${\cal J}: \chi(M) \times 
 C^{\infty}(M, \reals) \rightarrow \chi(M) \times C^{\infty}(M, \reals) $
 be a $C^{\infty}(M, \reals)$-linear map.  The Nijenhuis torsion of
  $ {\cal J}$ is defined as follows:

\begin{eqnarray*}
N_{ {\cal J}}\Big((X_1,f_1) , \ (X_2,f_2) \Big)& = [{\cal J}(X_1,f_1), \
 \ {\cal J}(X_2,f_2)]_{\star}
 -  {\cal J}[{\cal J}(X_1,f_1), \ (X_2,f_2) ]_{\star}\cr
&  -{\cal J}[(X_1,f_1), \ {\cal J}(X_2,f_2)]_{\star}
+ {\cal J}^2 [(X_1,f_1) , \ (X_2,f_2)]_{\star}.
\end{eqnarray*}

 \noindent  Then ${\cal J}$ is said to be a {\it recursion operator} of a 
 Jacobi structure $(\G,E)$ (see~\cite{mmp}) if 
 $$N_{ {\cal J}}\Big(\widetilde\#_{_{(\G,E)}}(\alpha_1,g_1), \widetilde\#_{_{(\G,E)}}(\alpha_2,g_2) \Big)=0$$ 

 \noindent for any $(\alpha_i,g_i) \in \Omega^1(M)\times C^{\infty}(M,\reals)$,
 $i=1,2$.

  \begin{dfn}~\cite{mmp}
 \label {mmp} Let $ {\cal J}: \chi(M) \times 
C^{\infty}(M, \reals) \rightarrow \chi(M) \times C^{\infty}(M, \reals) $
be a $C^{\infty}(M, \reals)$-linear map, and let $(\G, E)$ be a Jacobi structure on $M$. The triple $(\G,E,{\cal J})$ is said to be a Jacobi-Nijenhuis
structure on $M$ if ${\cal J}$ is a recursion operator of $(\G,E)$, if
  ${\cal J } \circ \widetilde \#_{_{(\G,E)}}
 = \widetilde \#_{_{(\G,E)}} \circ \ ^t  {\cal J}$ and if
  $(\G_1,E_1)$ is a Jacobi structure compatible with $(\G, E)$, 
 where $\G_1$ and $E_1$ 
 are characterized by 

$$\widetilde \#_{_{(\G_1,E_1)}}= {\cal J}\circ \widetilde \#_{_{(\G,E)}}.$$
\end{dfn}

In fact, a $C^{\infty}(M, \reals)$-linear map ${\cal J}$ from 
 $\chi(M) \times C^{\infty}(M, \reals)$ into itself is equivalent to 
 determining a quadruplet $(J, X_0, \alpha_0, \varphi_0)$, where $J$ is a
 $(1,1)$-tensor field on $M$, $X_0$ is a vector field, $\alpha_0$ is a 
 differential 1-form and $\varphi_0$ is a smooth function 
such that
 $$ {\cal J}(X,f)=(JX+fX_0, \ i_{X} \alpha_0+f \varphi_0).$$
  
  \noindent We set ${\cal J}=(J, X_0, \alpha_0, \varphi_0)$, the transpose
 $^t{\cal J}$  of ${\cal J}$ is defined by

  $$^t{\cal J}(\beta, g)= ( \ ^tJ \beta + g \alpha_0,\  i_{X_0} \beta
  + g \varphi_0), $$
 
\noindent for any $(\beta, g) \in \Omega^1(M) \times C^{\infty}(M, \reals)$.
  For any Jacobi structure $(\G,E)$, 

\[
 {\cal J } \circ \widetilde \#_{_{(\G,E)}}= \widetilde \#_{_{(\G,E)}} \circ \ ^t  {\cal J} \iff  
\left \{ \begin{array}{ll}
i_E \alpha_0=0\cr
JE=\G \alpha_0 + \varphi_0 E\cr
J\G\alpha -\G ^tJ\alpha= (i_E \alpha) X_0
 +(i_{X_0} \alpha) E, \ \forall \alpha \in \Omega^1(M) 
\end{array}  
\right. \]

\noindent From now on, we shall assume that $i_E \alpha_0=0$. \par

 We can now compare Jacobi-Nijenhuis structures with
  structures defined by a  $(1,1)$-tensor field $J$ and a Jacobi structure 
  which are compatible in the sense of Definition~\ref{Com}. \par

   We start with the simple example where 
   $\G=0$ (i.e. the Jacobi structure is 
 just a non-zero vector field on $M$), then for  any $(1,1)$-tensor field 
 $J$, the properties $(i)$, $(ii)$ and $(iii)$ of Definition~\ref{Com} 
 are satisfied. Taking into account the above 
 equivalence, we deduce that a recursion operator
  ${\cal J}=(J, X_0, \alpha_0, \varphi_0)$  of $(0,E)$ defines a
  Jacobi-Nijenhuis structure on $M$ if and only if $JE=\varphi_0 E$
 and  $X_0$ vanishes on the support of $E$. In other words, 
 any $(1,1)$-tensor field 
 $J$ is compatible with $(0,E)$ but for a recursion operator
  ${\cal J}=(J, X_0, \alpha_0, \varphi_0)$, the pair $(E, {\cal J})$ is a    
 Jacobi-Nijenhuis structure on $M$ if and only if the following conditions
   are fulfilled:

$$JE=\varphi_0 E \quad  {\rm and} \quad X_0=0 \ \mbox {\rm on supp}(E).$$

   We now move on to the general case. Let
 ${\cal J}=(J, X_0, \alpha_0, \varphi_0)$  be a recursion operator
 of a Jacobi structure $(\G,E)$,
  where $\G$ is not identically zero on $M$. There are two alternatives: \par
 
\noindent {\it First case:} $X_0=0$ on the support of $E$.
 We shall show that in such a case, $(\G, E, {\cal J })$ is a Jacobi-Nijenhuis 
 structure if and only if $JE=\G \alpha_0 + \varphi_0 E$ and $J$ fulfills
 the axioms $(i)$, $(ii)$ and $(iii)$ of Definition~\ref{Com}.\par 

 Assume that $(\G, E, {\cal J })$ is a Jacobi-Nijenhuis 
 structure. By hypothesis, we have $(i)$ and the equality
 $JE=\G \alpha_0 + \varphi_0 E$. Moreover, 
 for every integer $k \geq1$, the equation 
  $$\widetilde \#_{_{(\G_k,E_k)}}= {\cal J}^k\circ \widetilde \#_{_{(\G,E)}}$$
  has a unique solution that is $(\G_k, E_k)=(J^k \G, J^k E)$. 
  A result proven in~\cite{mmp} shows that if $(\G, E, {\cal J })$ is a 
 Jacobi-Nijenhuis structure then for any integer $k \geq1$, the pair
 $(J^k \G, J^k E)$ is a Jacobi-Nijenhuis structure that is compatible with
 $(\G, E)$. So, we obtain in particular $(iii)$.  
 Furthermore, setting 
  $${\cal V}=\{x \in M  | E(x)=0 \ \mbox{\rm in a neighborhood of} \ x \},$$
    it follows from the fact that ${\cal J}$ is a recursion
 operator that the Nijenhuis torsion $N_J$ of $J$ 
 satisfies the following equation on  ${\rm supp} (E)\cup {\cal V}$:

\begin{equation} \label{eq:Nij}
N_J(\G \alpha+fE,\G \beta+gE)=0.
\end{equation}  

 \noindent Since ${\rm supp} (E) \cup {\cal V}$ is a dense subset of $M$,
 by an argument of continuity, Equation (~\ref{eq:Nij}) is valid on the whole
 manifold $M$. Furthermore,  $JE=\G \alpha_0 + \varphi_0 E$ implies
 that $J^kE$ is also of the form $\G \alpha + \varphi E$, for all $k \geq 1$.
  So, 
 $$N_J(\G \alpha, J^k E)=0, \ \forall \alpha \in \Omega^1(M).$$
  
\noindent Applying Theorem~\ref{rec}, we obtain $(ii)$. 
  This shows that $J$ belongs to the particular class of structures satisfying 
  Definition~\ref{Com}  and $JE=\G \alpha_0 + \varphi_0 E$.
  
 Conversely, assume that  the relations $JE=\G \alpha_0 + \varphi_0 E$, 
  $(i)$, $(ii)$ and $(iii)$ are satisfied on $M$.
 Then  using Theorem~\ref{pencil} and Definition~\ref{mmp},
  we deduce that $(\G, E, {\cal J })$ is a 
 Jacobi-Nijenhuis structure.  \\

  \noindent {\it Second case:} if $X_0$ is not identically 
 zero on the support of $E$, then $J \circ \G-\G \circ ^tJ \ne 0$.
  However,  $(\G, E, {\cal J })$ may be a Jacobi-Nijenhuis structure.\\

 We summarize our discussion in the following proposition:

 \begin{prop}Let $(M,\G, E)$ be a Jacobi manifold and let $J$ be a 
 (1,1)-tensor field on $M$ such that $J \circ \G=\G \circ \ ^tJ$.
 Assume that there exist a vector field $X_0$,
 a 1-form $\alpha_0$ and a smooth function $\varphi_0$ such that
  ${\cal J}=(J, X_0, \alpha_0, \varphi_0)$ 
 is a recursion operator of $(\G,E)$. Then $(\G, E, {\cal J })$ is a 
 Jacobi-Nijenhuis structure if and only if the following conditions
 are satisfied:

 \begin{itemize}
  \item[$(C_1)$] \ \ \ \ \ \ $JE=\G \alpha_0 + \varphi_0 E$,
 \item[$(C_2)$] \ \ \ \ \ \   $ [J^kE, \G ]+[E, J^k\G]=0,$ 
 for any integer $k \geq 1$,
\item[$(C_3)$] \ \ \ \ \ \  $\big \langle \ ^t J \gamma, \ \G 
 \Big(C(\G,J)(\alpha, \beta) \Big) - \G(\alpha, \beta)JE 
 +\G(\alpha, \ ^t J\beta)E \big \rangle=0,$
 for any $\alpha, \ \beta, \ \gamma \in \Omega^1(M)$.
\end{itemize}
\end{prop}

\section { Nijenhuis tensors and homogeneous Poisson structures}

\begin{dfn}
A homogeneous Poisson manifold $(M,\pi, Z)$ is a Poisson manifold $(M,\pi)$ with a vector field $Z$ over $M$ such that 

$$[Z,\pi]=-\pi.$$
\end{dfn}

\begin{thm}
Assume that $(M,\pi, Z)$ is a homogeneous Poisson manifold. Let $J$ be 
Nijenhuis tensor compatible with $\pi$. Then $(M, J\pi, Z)$ is a homogeneous 
Poisson manifold if and only if the following property is satisfied:

\begin{equation} \label{eq:comp}
\pi \circ (L_Z \circ \ ^tJ- ^tJ \circ L_Z)=0,
\end{equation}

\noindent  where $L_Z=i_Zd+di_Z$ is the Lie derivation by $Z$. When this property holds,  $J\pi -\lambda \pi$ defines a Poisson 
pencil which is homogeneous with respect to $Z$.
\end{thm}

\noindent {\it Proof:} Taking into account 
Theorem~\ref{pencil},  we have only to prove that $[Z,J\pi]=-J\pi$. Let us compute $[Z,J\pi]$. We obtain

\begin {eqnarray*}
[Z,J\pi](df,dg) &=& L_Z \big(J\pi(df,dg)\big) -J\pi\big(L_Zdf,dg\big) -
J\pi\big(df,L_Zdg\big) \\
 &=& L_Z \big(\pi(^tJdf,dg) \big) -\pi\big(^tJL_Zdf,dg\big) -
\pi\big(^tJdf,L_Zdg\big) 
\end{eqnarray*}

\noindent Since 
$$L_Z \big(\pi(^tJdf,dg) \big)=[Z,\pi]\big(^tJdf,dg\big) +
\pi\big(L_Z \ ^tJdf,dg\big) + \pi\big(^tJdf, L_Zdg\big), $$

\noindent we obtain 
\begin {eqnarray*}
[Z,J\pi](df,dg) &=&[Z,\pi]\big(^tJdf,dg\big) +\pi\big(L_Z \ ^tJdf,dg\big)
-J\pi\big(L_Zdf,dg\big)\\
&=&-\pi\big(^tJdf,dg\big) +\pi\big(L_Z \ ^tJdf,dg\big)-J\pi\big(L_Zdf,dg\big)
\end{eqnarray*}

\noindent Hence, the relation $[Z,J\pi]=-J\pi$ is equivalent to the following 
one:

$$\pi \circ L_Z \circ \ ^tJ= \pi \circ \ ^tJ \circ L_Z.$$
\noindent This proves the theorem. \hfill \qed

\begin{dfn} A homogeneous Poisson manifold $(M,\pi, Z)$ equipped with a
 Nijenhuis tensor $J$  which is compatible with $\pi$ and satisfies
 equation (~\ref{eq:comp}) is said to be a homogeneous Poisson-Nijenhuis
 manifold.
\end{dfn}

\begin{cor}Let $(M,\pi, J)$ be a Poisson-Nijenhuis manifold. If $\pi$ is homogeneous with respect to a vector field $Z$ and  if the  following property holds
\begin{equation}
[Z,JX]=J[Z,X], \quad  \quad \forall X \in \chi(M), 
\end{equation}
\noindent then the triple  $( M, \pi, J)$ is a homogeneous Poisson-Nijenhuis 
 manifold with respect to $Z$.
\end{cor}

\noindent {\it Proof:}  We obtain this corollary using the above theorem and 
the fact that
$$
[Z,JX]=J[Z,X], \quad \forall X \in \chi(M) \iff \quad  L_Z \circ ^tJ=^tJ \circ L_Z.$$

\begin{dfn} A map $\psi: (M_1,\G_1,E_1) \rightarrow (M_2,\G_2,E_2)$
 between two Jacobi manifolds is said to be a {\it conformal Jacobi morphism}
if  there exists a function $a \in C^{\infty}(M_1, \reals)$ which vanishes nowhere such that for any $f,g \in C^{\infty}(M_2, \reals)$ we have:
$$\{a (f \circ \psi), a (g \circ \psi) \}_1= a (\{f,g\}_2 \circ \psi),$$
\noindent where the brackets $\{ \ , \ \}_1$ and $\{ \ , \ \}_2$ are the Jacobi
brackets associated with  $(\G_1,E_1)$ and $(\G_2,E_2)$ respectively.
\end{dfn}

Homogeneous Poisson manifolds are closely related to Jacobi manifolds and
 their relations were established in [D-L-M]. In terms of Poisson pencils,
we have the following results.

\begin{prop}
\label{hom-jac}
Let  $\{ .,.\}_{_{\lambda}}$ be a Jacobi pencil on $M$. There exists a 
Poisson pencil on $M \times \reals$ such that the projection
$P: M \times \reals \rightarrow M$ is a conformal Jacobi morphism, for each 
 $\lambda$.
\end{prop}

\noindent {\it Proof:} If $(\G_i,E_i)$ denotes the Jacobi structure on $M$ associated to $\{.,. \}_{_i}$, with $i=1,2$, then according to 
 Proposition~\ref{car}, the Poisson pencil on  $M \times \reals$ is given by
$\pi_1 - \lambda \pi_2$ where

$$\pi_i(x,t)=e^{-t}\Big(\G_i(x)+{\partial \over \partial t} \wedge E_i \Big).$$
\noindent One may easily verify that $P:(M \times \reals , \pi_{\lambda})
 \rightarrow (M,\{.,.\}_{_{\lambda}}) $ is a conformal Jacobi morphism.
\hfill \qed

 Conversely, we may prove that homogeneous Poisson pencils give Jacobi pencils by using a proof from~\cite{dlm}. More precisely, we have:
\begin{prop}
\label{hom-pen}
Let  $\pi_{\lambda}$ be a homogeneous Poisson pencil on $M$ with 
respect to the vector field $Z$, and let $N$ be a sub-manifold of $M$ of 
codimension 1 which is transverse to $Z$. Then there exists a 
Jacobi pencil on $N$ such that for any pair of functions 
$(f,g)$ defined on an open set $U$ of $M$, satisfying $<Z, df>=f$ and
$<Z,dg>=g$, we have
$$\{f_{|_{N\cap U}},g_{|_{N\cap U}}\}_{\lambda}=\pi_{\lambda}(df, dg)_{|_{N\cap U}}.$$
\end{prop}

\begin{cor}
\label{ed}
Let $(M, \G, E)$ be a Jacobi manifold and let $J$ be a Nijenhuis tensor on 
$M$, which is compatible with $(\G,E)$. Then there exists a sequence of
 Poisson-Nijenhuis
structures $(\pi_k)$ on $M \times \reals$ such that the projection
$P_k: (M \times \reals, \pi_k) \rightarrow (M, \G, E)$ is a conformal Jacobi morphism,  for each $k \geq 1$ . \par

Conversely, if $(M, \pi, J)$ is a homogeneous Poisson-Nijenhuis manifold with respect to the vector field $Z$ and if $N$ is a sub-manifold of $M$ of codimension 1, which is transverse to $Z$, then there exists a sequence of 
pairwise compatible Jacobi structures on $N$ determined by $\pi$, $Z$ and $J$.

\end{cor}

\noindent This corollary is a direct consequence of  Theorem~\ref{pencil}
  and Propositions~\ref{hom-jac} and~\ref{hom-pen}.

\bigskip

\noindent {\bf Acknowledgments.} I express my gratitude to the 
Abdus Salam International Centre for Theoretical Physics for its support. 
 I wish to thank R. Brouzet, J.-P. Dufour and A. Kuku for helpful 
 discussions. Thanks are due also to C.-M. Marle for pointing out 
 Ref.~\cite{mmp}, J. Monterde and J. M. Nunes da Costa for providing me
 with Ref.~\cite{mmp} and Ref.~\cite{nunes}, respectively. 
  I am indebted to the referee for useful comments
 which enabled me to bring this paper to its present form.

\bigskip

\noindent {\bf Current address:} \par
\noindent {\small DEPARTMENT OF MATHEMATICS,} 
 {\small UNIVERSITY OF NORTH CAROLINA,\par
  \noindent CHAPEL HILL, NC} 27599-3250.\par
  \noindent E-mail: aissaw@math.unc.edu

\end{document}